\documentclass[Journal]{IEEEtran}
\usepackage{cite}
\usepackage{epstopdf}
\usepackage{subfigure}
\usepackage{bm}
\usepackage{amsmath}
\usepackage{amssymb}
\usepackage{pifont}
\usepackage{xtab}
\usepackage{mathrsfs}
\usepackage{subeqnarray}
\usepackage{cases}
\usepackage{algorithm}
\usepackage{algcompatible}

\algblockdefx{FORALLP}{ENDFAP}[1]%
  {\textbf{for all }#1 \textbf{do in parallel}}%
  {\textbf{end for}}

\ifCLASSINFOpdf

\else

   \usepackage[dvips]{graphicx}
  \graphicspath{{../eps/}}
   \DeclareGraphicsExtensions{.eps}
\fi

\IEEEoverridecommandlockouts

\newcommand{\abs}[1]{\lvert#1\rvert}

\newtheorem {Condition1} {\bf{Cond.}}
\newtheorem {Equivalence1} {\bf{Lemma}}
\newtheorem {Assumption1} {\bf{Asmp.}}

\begin{document}

\title{Improved Sufficient Conditions for Exact Convex Relaxation of Storage-Concerned ED}

\author{
\IEEEauthorblockN{Chao Duan,~\IEEEmembership{Student Member,~IEEE}, Lin Jiang,~\IEEEmembership{Member,~IEEE}, Wanliang Fang, Xin Wen, Jun Liu,~\IEEEmembership{Member,~IEEE}}

\thanks{This work was supported in part by Engineering and Physical Sciences Research Council (EP/L014351/1).}

\thanks{C. Duan, W. Fang and J. Liu are with the Department of Electrical Engineering, Xi¡¯an Jiaotong University, Xi¡¯an 710049, China.  C. Duan is also with Department of Electrical Engineering and Electronics, University of Liverpool, Liverpool L69 3GJ, U.K. }

\thanks{ L. Jiang and X. Wen are with Department of Electrical Engineering and Electronics, University of Liverpool, Liverpool L69 3GJ, U.K.}

}

\maketitle

\begin{abstract}
To avoid simultaneous charging and discharging of storages, complementarity constraints are introduced to storage-concerned economic dispatch (ED), which makes the problem non-convex. This letter concerns the conditions under which the convex relaxation of storage-concerned ED with complementarity constraints is exact. Two new sufficient conditions are proposed, proved and verified to significantly reduce the conservatism of recent results \cite{lisufficient,li2015further}.
\end{abstract}


\IEEEpeerreviewmaketitle

\vspace{-0.12cm}
\section{Introduction}
One difficulty appears in storage-concerned economic dispatch (ED) calculation is how to avoid simultaneous charging and discharging, which is unrealistic for most energy storage technologies \cite{jabr2015robust}. The introduction of auxiliary binary variables \cite{malysz2014optimal,jabr2015robust} as well as complementarity constraints \cite{lisufficient} are widely considered in the literature. The former leads to mixed integer programming (MIP) and the latter results in non-convex nonlinear programming (NLP). Both are NP-hard. By simply dropping the complementarity constraints in the latter approach, the non-convex problem can be relaxed to a convex problem which is polynomial time solvable by interior point method. This letter concerns the conditions under which above relaxation is exact in the sense that the convex problem attains the same optimal solution as the original non-convex one. We improve the results in recent papers \cite{lisufficient,li2015further}. First, we propose a local marginal price (LMP) related sufficient condition which is weaker than those given in \cite{lisufficient} and \cite{li2015further}. Second, we present an even weaker condition concerning the sizes of the storages where and when the first condition is violated.

\vspace{-0.08cm}
\section{Storage-Concerned ED}
Consider a power network with bus set $\mathcal{N}$ and branch set $\mathcal{L}$. $\mathcal{T}=\{1,\dots,T \}$ denotes the set of time slots. The storage-concerned ED problem minimize the objective function
\begin{equation}\label{OPF_Obj}
v(\varOmega)= \sum_{t \in \mathcal{T}} \sum_{i \in \mathcal{N}}  \big(g_{i}(p_{i}^{d}(t))-f_{i}(p_{i}^{c}(t))+h_{i}(p_{i}^{g}(t))\big)
\end{equation}
where $\varOmega=\left( p_{i}^{c}(t),p_{i}^{d}(t),p_{i}^{g}(t) \right)_{\forall t \in \mathcal{T},i \in \mathcal{N}}$, subject to:
\begin{align}
&0\leq p_{i}^{c}(t) \leq \overline{P}_{i}^{c}, \quad \alpha_{i,1}(t),\alpha_{i,2}(t)\label{ESScharge}\\
&0\leq p_{i}^{d}(t) \leq \overline{P}_{i}^{d}, \quad \alpha_{i,3}(t),\alpha_{i,4}(t) \label{ESSdischarge}
\end{align}
\begin{align}
&p_{i}^{d}(t)p_{i}^{c}(t)=0\label{nosimcharge}
\\
&\underline{P}_{i}^{g} \leq p_{i}^{g}(t) \leq \overline{P}_{i}^{g}\label{OPFconstr_genp}\\
&\underline{S}_{i} \leq s_{i}(t) \leq \overline{S}_{i}, \quad \beta_{i,1}(t),\beta_{i,2}(t)\label{OPFconstr_ESSstore}\\
&R_{i}^{d}\Delta t \leq p_{i}^{g}(t+1)-p_{i}^{g}(t) \leq R_{i}^{u}\Delta t \label{Ramprate}\\
&\sum_{i\in \mathcal{N}} \left( p_{i}^{g}(t)+p_{i}^{d}(t)-p_{i}^{c}(t) \right)=\sum_{i\in \mathcal{N}} D_{i}(t), \quad \lambda(t)\label{Powerbalance}\\
\begin{split}\label{Linelimit}
\underline{P}_{j}^{l} \leq \sum_{i\in \mathcal{N}}GSF_{j-i}\left( p_{i}^{g}(t)+p_{i}^{d}(t)-p_{i}^{c}(t)- D_{i}(t) \right)\\
\leq \overline{P}_{j}^{l}, \quad \mu_{j,1}(t),\mu_{j,2}(t)
\end{split}
\end{align}
where
\begin{equation}\label{StorgedEnergy}
\begin{split}
s_{i}(t)=&(1-\varepsilon_{i})^{t}s_{i}^{0}\\
&+\sum_{\tau=1}^{t}(1-\varepsilon_{i})^{t-\tau}\left( \eta_{i}^{c}p_{i}^{c}(\tau)-p_{i}^{d}(\tau)/\eta_{i}^{d} \right)\Delta t.
\end{split}
\end{equation}
The decision variables include the grid-side energy storage charging power $p_{i}^{c}(t)$ and discharging power $p_{i}^{d}(t)$ and the generator active power output $p_{i}^{g}(t)$. Convex quadratic discharging cost $g_{i}$, linear storage charging fee $f_{i}$ and convex quadratic generation cost $h_{i}$ form the objective function (\ref{OPF_Obj}). Inequalities (\ref{ESScharge}) and (\ref{ESSdischarge}) set the limits for storage charging and discharging power. Complementarity constraint (\ref{nosimcharge}) ensures storages operate either in the charge or discharge mode. The upper and lower bounds of generator output and storage energy are enforced by (\ref{OPFconstr_genp}) and (\ref{OPFconstr_ESSstore}). (\ref{Ramprate}) represents the generator ramp rate constraint. (\ref{Powerbalance}) is the power balance equation of the whole system, and (\ref{Linelimit}) represents the bidirectional transmission capacity limits. The self-discharging effect has been considered by the self-discharge rate $\varepsilon_{i}$. $\alpha_{i,1}(t)$, $\alpha_{i,2}(t)$, $\alpha_{i,3}(t)$, $\alpha_{i,4}(t)$, $\beta_{i,1}(t)$, $\beta_{i,2}(t)$, $\lambda(t)$, $\mu_{j,1}(t)$ and $\mu_{j,2}(t)$ are multipliers of corresponding constraints. For brevity, we do not include net charging requirements constraints and time-varying limits \cite{lisufficient} in our formulation. But the propositions in this paper also hold when those constraints are considered.

\section{Convex Relaxation and Exactness}
Two problems are considered. The first problem is the non-convex original problem formally stated as OP: $\underset {\Omega}{\mathrm{min}}\ v(\Omega)$ s.t. (2)$\thicksim$(9) whose feasible set, optimum and optimal solution are denoted as $\mathcal{F}_{0}$, $v^{*}_{0}$ and $\varOmega^{*}_{0}$. We assume OP has unique global solution. The second problem is the convex relaxation problem RP: $\underset {\Omega}{\mathrm{min}}\ v(\Omega)$ s.t. (2) $\thicksim$ (3), (5) $\thicksim$ (9) with feasible set, optimum and optimal solution denoted as $\mathcal{F}_{1}$, $v^{*}_{1}$ and $\varOmega^{*}_{1}$. Obviously, $\mathcal{F}_{0}\subseteq\mathcal{F}_{1}$ and $v^{*}_{1} \leq v^{*}_{0}$.

At first, we present an improved LMP related condition and the exactness of RP under this condition.
\begin{Condition1}
$\forall i \in \mathcal{N}$, $\forall t \in \mathcal{T}$,
\begin{equation}\label{LMPbound}
LMP_{i}(t)> \frac{f'_{i}(p_{i}^{c}(t))-\eta_{i}^{c}\eta_{i}^{d}g'_{i}(p_{i}^{d}(t))}{1-\eta_{i}^{c}\eta_{i}^{d}}
\end{equation}
where $LMP_{i}(t)=\lambda(t)+\sum_{j\in \mathcal{L}}GSF_{j-i}\left( \mu_{j,1}(t)-\mu_{j,2}(t) \right)$.
\end{Condition1}

\begin{Equivalence1}
If Cond. 1 holds for the primal and dual solution of RP, $v^{*}_{1} = v^{*}_{0}$ and $\varOmega^{*}_{1} = \varOmega^{*}_{0}$.
\end{Equivalence1}
\begin{proof}
We already have $\varOmega^{*}_{0} \in  \mathcal{F}_{1}$ and $v^{*}_{1} \leq v^{*}_{0}$. It suffices to show $\varOmega^{*}_{1} \in  \mathcal{F}_{0}$ and $v^{*}_{0} \leq v^{*}_{1}$. Assume that $\exists$ $p_{i}^{c}(t)>0$ and $p_{i}^{d}(t)>0$ for some $i \in \mathcal{N}$, $t \in \mathcal{T}$ in $\Omega^{*}_{1}$. Considering the KKT optimality condition of RP, we have $\alpha_{i,1}(t)=0$, $\alpha_{i,3}(t)=0$, $\alpha_{i,2}(t)\geq0$ and $\alpha_{i,4}(t)\geq0$. In addition, $\frac{\partial L}{\partial p_{i}^{c}(t)}=\frac{\partial L}{\partial p_{i}^{d}(t)}=0$, i.e.
\begin{equation}\label{KKT1}
-f'_{i}(p_{i}^{c}(t))+\alpha_{i,2}(t)-\eta_{i}^{c}\varGamma(t) \Delta t+LMP_{i}(t)=0
\end{equation}
\begin{equation}\label{KKT2}
g'_{i}(p_{i}^{d}(t))+\alpha_{i,4}(t)+\varGamma(t) \Delta t / \eta_{i}^{d}-LMP_{i}(t)=0
\end{equation}
where $\varGamma(t)=\sum_{\tau=t}^{T}(1-\varepsilon_{i})^{(\tau-t)}\left( \beta_{i,1}(\tau)-\beta_{i,2}(\tau) \right)$. Eliminating $\varGamma(t) \Delta t$ from (\ref{KKT1}) and (\ref{KKT2}), we obtain
\begin{equation}
LMP_{i}(t)=\frac{f'_{i}(p_{i}^{c}(t))-\eta_{i}^{c}\eta_{i}^{d}g'_{i}(p_{i}^{d}(t))}{1-\eta_{i}^{c}\eta_{i}^{d}}-\frac{\alpha_{i,2}(t)+\eta_{i}^{c}\eta_{i}^{d}\alpha_{i,4}(t)}{1-\eta_{i}^{c}\eta_{i}^{d}}
\end{equation}
which contradicts to Cond. 1. So $p_{i}^{c}(t)p_{i}^{d}(t)=0$, $\forall i \in \mathcal{N}$, $t \in \mathcal{T}$, i.e. $\varOmega^{*}_{1} \in  \mathcal{F}_{0}$ and $v^{*}_{0} \leq v^{*}_{1}$. Therefore $v^{*}_{0} = v^{*}_{1}$. Due to the uniqueness of global solution, we have $\varOmega^{*}_{1} = \varOmega^{*}_{0}$.
\end{proof}

It is easy to verify that Cond. 1 is strictly weaker than the conditions proposed in \cite{lisufficient,li2015further}. In particular, if $f'_{i}(p_{i}^{c}(t))<\eta_{i}^{c}\eta_{i}^{d}g'_{i}(p_{i}^{d}(t))$, RP is exact even when the LMP is negative (a well-known situation for simultaneous charging and discharging \cite{jabr2015robust}) which is not allowed in the conditions in \cite{lisufficient,li2015further} with positive $f'_{i}(p_{i}^{c}(t))$. Since the LMP can be predicted based on historical data, cond. 1 can be checked before solving RP.

To establish exactness under a weaker condition, we need to make the following assumption.
\begin{Assumption1}
$f'_{i}(p_{i}^{c}(t))< g'_{i}(p_{i}^{d}(t))$, $\forall i \in \mathcal{N}$, $t \in \mathcal{T}$.
\end{Assumption1}

If the storages are owned by the grid, $f'_{i}(p_{i}^{c}(t))<0$ and $g'_{i}(p_{i}^{d}(t))>0$. Asmp. 1 holds trivially. If the storages are not own by the grid, Asmp. 1 means storage charging prices paid to the grid are universally less than the discharging compensation prices paid to the storage owners, which is the prerequisite for the storage owners to participate in economic dispatch. Asmp. 1 is also used in \cite{lisufficient} and \cite{li2015further}.

Then we give a weaker exactness condition concerning the the sizes of the storages installed at buses with low LMPs.
\begin{Condition1}
$\forall i \in \mathcal{N}$, $ t \in \mathcal{T}$ at which the dual solution of RP violates Cond. 1, $s_{i}(\tau)<\overline{S}_{i}$, $\forall \tau \geq t$.
\end{Condition1}

\begin{Equivalence1}
Under Asmp. 1, if Cond. 2 holds for the primal and dual solution of RP, $v^{*}_{1} = v^{*}_{0}$ and $\varOmega^{*}_{1} = \varOmega^{*}_{0}$.
\end{Equivalence1}

\begin{proof}
If $\nexists$ $ i \in \mathcal{N}$, $ t \in \mathcal{T}$ at which the dual solution of RP violates (\ref{LMPbound}), Cond. 1 holds, the exactness is proved by Lemma 1. If $\exists$ such $i \in \mathcal{N}$, $ t \in \mathcal{T}$, we assume that $p_{i}^{c}(t)>0$ and $p_{i}^{d}(t)>0$ in $\Omega^{*}_{1}$. Hence (\ref{KKT1}) and (\ref{KKT2}) hold with $\alpha_{i,2}(t)\geq0$ and $\alpha_{i,4}(t)\geq0$. Eliminating $LMP_{i}$ from (\ref{KKT1}) and (\ref{KKT2}), we obtain
\begin{equation}
-f'_{i}(p_{i}^{c}(t))+g'_{i}(p_{i}^{d}(t))+\alpha_{i,2}(t)+\alpha_{i,4}(t)+(\frac{1}{\eta_{i}^{d}}-\eta_{i}^{c})\varGamma(t) \Delta t=0
\end{equation}
Considering Asmp. 1 and noticing the positivity of $\alpha_{i,2}(t)$, $\alpha_{i,4}(t)$ and $1/\eta_{i}^{d}-\eta_{i}^{c}$, we have $\varGamma(t)<0$. So $\exists \ \tau \geq t$, $\beta_{i,2}(\tau)>0$, i.e. $s_{i}(\tau)=\overline{S}_{i}$. This contradicts to Cond. 2. So $p_{i}^{c}(t)p_{i}^{d}(t)=0$, $\forall i \in \mathcal{N}$, $t \in \mathcal{T}$, i.e. $\varOmega^{*}_{1} \in  \mathcal{F}_{0}$ and $v^{*}_{0} \leq v^{*}_{1}$. Therefore $v^{*}_{0} = v^{*}_{1}$. Similarly, we have $\varOmega^{*}_{1} = \varOmega^{*}_{0}$.
\end{proof}

Lemma 2 states that if the LMPs at some buses go below the bound given in Cond. 1, the exactness of RP can still be guaranteed provided the storage installed at those buses has large enough energy capacity. The maximal energy capacity needed can also be estimated based on the forecasted LMPs. It can be analyzed from (\ref{KKT1}) and (\ref{KKT2}) that, by assuming $\varGamma(t) \geq 0$, the storage charges with $\eta_{i}^{c}\overline{P}_{i}^{c}$ when $LMP_{i}<f'_{i}$, and discharge with $\text{min} \left(\overline{P}_{i}^{d},\frac{(1-\varepsilon_{i})s_{i}(t-1)-\underline{S}_{i}}{\Delta t} \right)/\eta_{i}^{d}$ when $LMP_{i}>g'_{i}$.

The following procedure can be used to check the exactness of RP. 1) Forecast the LMPs by using historical data and check Cond.1; if satisfied, RP is exact; 2) If Cond.1 does not hold, find out those $i$ and $t$ at which Cond.1 is violated; 3) Estimate the maximum energy capacity needed at bus $i$. If $\overline{S}_{i}$ is larger than the maximal energy capacity needed, RP is still exact.

\section{Numerical Validation}
Numerical tests are conducted on IEEE 30-bus system with 3 wind farms and 5 storages. The maximal load is 189 MW, and the maximal wind generation is 60MW. Both load and wind generation vary according to daily forecasted curves. A time horizon of 24h in time steps of 0.5h is considered. The relaxed problem is solved by SDPT3 with YALMIP.

Results are shown in Table I. When fix $(f'_{i},g'_{i})=(1.5,2.5)$, the lower bound of LMP for the exactness of RP is 1.5 provided by \cite{lisufficient} versus -2.76 by Cond. 1 in this paper. Row 1 to row 3 of Table I demonstrate the validity and superiority of Cond. 1. When LMPs decrease to -3, Cond.1 and Cond.2 are violated thus simultaneous charging and discharging happens, shown in row 4. But after we enlarge the energy capacity of storages from 2 MWh to 10 MWh, the exactness of RP is recovered, shown in row 5. The proposed Cond. 1 and Cond. 2 significantly reduce the conservatism of previous results.

\begin{table}\centering
\caption{Numerical Validation of Exactness Conditions}\label{tab1}
\begin{minipage}[c]{250pt}
\begin{tabular}
{|c|c|c|c|c|c|} \hline
$(f'_{i},g'_{i})$ & $LMP_{i}$ & Cond.1 & Cond.2 & Conds in \cite{lisufficient} & max($\abs{p^{c}_{i}p^{d}_{i}}$) \\ \hline
(1.5,2.5) & 2.0 & yes & $\diagup$ & yes & 1.4e-20   \\ \hline
(1.5,2.5) & 1.2 & yes & $\diagup$ & no & 3.8e-11   \\ \hline
(1.5,2.5) & -1 & yes & $\diagup$ & no & 8.3e-11   \\ \hline
(1.5,2.5) & -3 & no & no & no & \bf{0.0012}   \\ \hline
(1.5,2.5) & -3 & no & yes & no &  1.4e-12  \\ \hline
\end{tabular}
\vspace{-9pt}\renewcommand{\footnoterule}{}
\end{minipage}
\end{table}

\appendices


\bibliographystyle{ieeetr}

\begin{thebibliography}{100}

\bibitem{jabr2015robust}
R.~Jabr, S.~Karaki, and J.~Korbane, ``Robust multi-period OPF with storage and
  renewables,'' \emph{IEEE Trans. Power Syst.}, vol.~PP, no.~99,
  pp.~1--10, 2014.

\bibitem{malysz2014optimal}
P.~Malysz, S.~Sirouspour, and A.~Emadi, ``An optimal energy storage control
  strategy for grid-connected microgrids,'' \emph{IEEE Trans. Power Syst.}, vol.~5, pp.~1785--1796, July 2014.

\bibitem{lisufficient}
Z.~Li, Q.~Guo, H.~Sun, and J.~Wang, ``Sufficient conditions for exact
  relaxation of complementarity constraints for storage-concerned economic
  dispatch,'' \emph{IEEE Trans. Power Syst.}, vol.~PP, no.~99,
  pp.~1--2, 2015.

\bibitem{li2015further}
Z.~Li, Q.~Guo, H.~Sun, and J.~Wang, ``Further discussions on sufficient
  conditions for exact relaxation of complementarity constraints for
  storage-concerned economic dispatch,'' {\em arXiv preprint arXiv:1505.02493},
  2015.


\end{thebibliography}

\end{document}